\documentclass[letterpaper, 10 pt, conference]{ieeeconf}  

\IEEEoverridecommandlockouts                              

\usepackage{graphicx} 

\title{\LARGE \bf
Optimization of 3-D flight trajectory of variable trim kites \\ for airborne wind energy production
}

\author{Rafal Noga$^{1}$, Xaver Paulig$^{1}$, Lukas Schmidt$^{1}$, Benjamin Karg*$^{1}$, Manfred Quack$^{1}$ and Mahmoud Soliman$^{1}$
\thanks{*corresponding author, 
	{\tt\small benjamin.karg@skysails.de}}%
\thanks{$^{1}$SkySails Power GmbH, Wendenstrasse 375,
	20537 Hamburg, Germany}%
\thanks{A pilot system of the SKS PN-14 is part of the "SkyPower100" project funded by the German Federal Ministry of Economics and Energy.}%
}

\begin{document}

\maketitle
\thispagestyle{empty}
\pagestyle{empty}

%
%
%
%
\section{SkySails Power PN-14}
SkySails Power GmbH is the leading manufacturer of light and efficient power kites that harness the wind's untapped supplies at high altitudes, aiming at profoundly altering wind energy's impact in achieving the global energy transition \cite{www_skysails}.

SkySails Power has developed a revolutionary onshore Airborne Wind Energy System, the {PN-14}, that is compact, light, and easy to manufacture, transport and operate. 
Subject to the site specific configuration, the PN-14 design enables rated power output up to 200 kW, using relatively large kites having area of 90-180~m2.  

\begin{figure}[h]
 	\includegraphics[width=\linewidth,trim={3cm 0 3cm 0},clip]{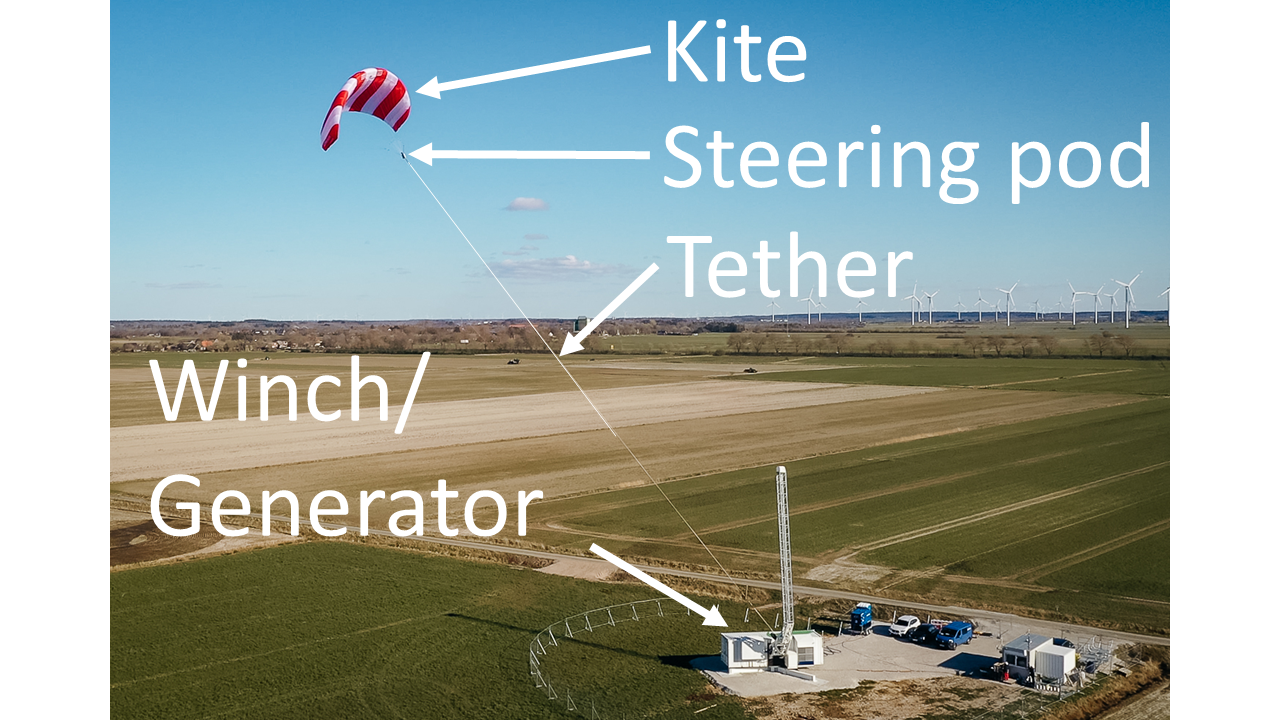}
 	\caption{Skysails PN-14}
 	\label{fig:pn14}
 \end{figure}

PN-14's main components are the ground station, built around a 30 ft container and housing a winch/generator, the 14~mm diameter tether connecting it to the steering pod that in turn is connected to the kite through a bridle of lines, see Fig. \ref{fig:pn14}. 
The principle of operation is simple: the tether is reeled out from the ground station under strong pulling force, thus generating significant amount of work $W_{reel\ out}$. 
After reaching the maximal length, it is reeled back in at lower force levels and thus at lower work expense $W_{reel\ in}$. 
Thus the mean net system power output $P_{output, mean} = (W_{reel\ out} - W_{reel\ in})/T_{cycle}$, with the duration of the power generation cycle  $T_{cycle}$.

The aerodynamic force on a tethered kite may be huge as it may fly few times faster than wind speed \cite{lloyd1980}. The exact value of the force depends on a number of parameters like wing aerodynamic coefficients, wind speed, kite position with respect to the wind direction and the rate of change of the tether length.
That is why in order to maximize mean cycle power, the kite flies a sophisticated trajectory across the sky, passing through extreme values of the parameters.

In order to further increase the power output, novel, variable trim kites have been developed for the PN-14. 
Varying the trim, corresponding to the angle of attack of the kite airfoil, allows to significantly modulate the aerodynamic coefficients of the kite, enabling operation within a larger range of forces pulling on the tether and thus resulting in significantly improved overall efficiency of the system \cite{fimpel23}.

However, the operation of variable trim kites differs significantly from that of previous constant trim kite generations.
It is also more complex and thus its mastering is one of the keys to a successful operation of the novel system \cite{soliman24}.

\section{Optimization of the kite trajectories}
Due to the very nature of the power generation cycle, tethered kites are operated in a wide range of frequently changing conditions and configurations that strongly alter their behavior, as described in the previous section.
Moreover, in order to maximize the Levelized Cost of Electricity, the system has to operate in a wide range of wind speeds and wind turbulence levels.
The main goal of the system operation is to maximize the electrical energy production in the whole range of the operating conditions, while satisfying a large number of operational constraints.
Some of the constraints are geometrical, related to shape of authorized airspace or presence of obstacles, some represent performance limitations of various feedback loops, other correspond to operational, aerodynamic, mechanical and electrical limits of various system components. 
Manipulated variables are the tether length, kite canopy direction and its trim.
Measured variables are physical quantities related to kite position relative to the ground station and to the power generation process. 
Ever changing wind speed spatial distribution represents a strong unmeasured disturbance to the kite traversing large volumes of space and it has to be estimated.
Thus optimizing the system performance in the whole range of operational conditions is a non-trivial task.  

Numerical optimization has been applied in the past to study the complex problem of optimization of system operation of tethered kites operated in a power generating cycle at SkySails and elsewhere \cite{zamm15} \cite{lago16} \cite{berra21}. 
However, kites considered in that work had constant aerodynamic coefficients.

Building on that work, we have developed a numerical optimization of the operation of the PN-14 using variable trim kites, taking advantage of variable aerodynamic coefficients.

The optimization problem is formulated as a classic optimal control problem with an economic objective.
We addressed all relevant constraints including the non-linear system dynamics. Most of the constraints  were implemented directly, some were approximated in a form more suitable for optimization and relaxation techniques were applied where needed.
Moreover, there was no need to add any significant artificial constraints to improve the performance of the numerical solution process. 

The optimization problem has been implemented in Python using CasADi for modeling using symbolic mathematics, algorithmic differentiation and interfacing to the non-linear programming solver IPOPT \cite{casadi}. IPOPT implements the interior point method \cite{ipopt}. 
An important aspect of the implementation is a simple software architecture with clean interfaces between software components, making it accessible to non-expert users, e.g. system engineers. 

The iterative solver is initialized using a realistic, synthetic, wind speed dependent solution guess.  

In order to cope with the non-convex nature of the non-linear optimization problem having a rather large number of complex non-convex constraints, globalization techniques have been used to reduce the probability of returning a local optimum significantly worse than the global one.

\section{Numerical results and discussion}
\begin{figure}
	\includegraphics[width=\linewidth,trim={17cm 3cm 17cm 16cm},clip]{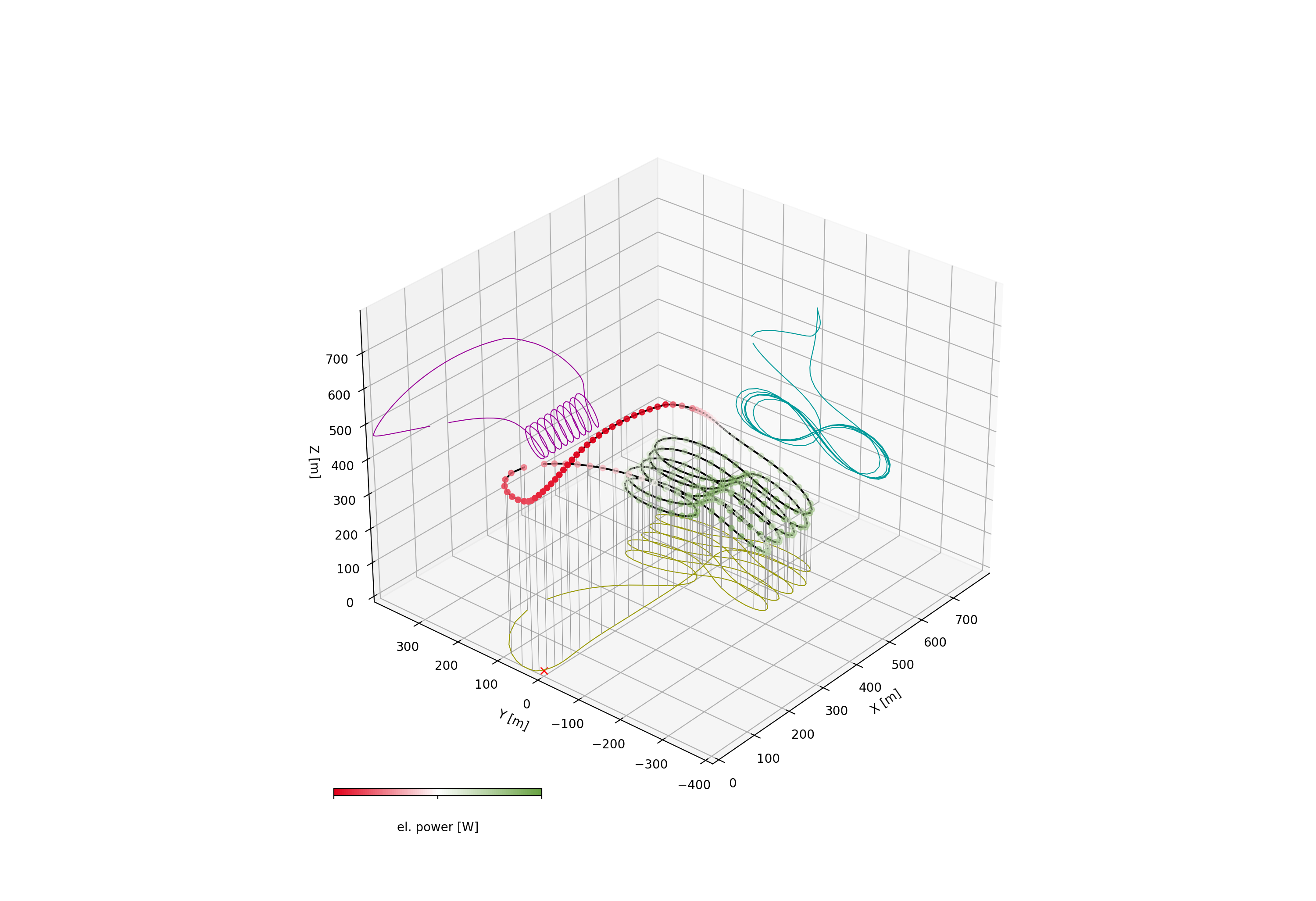}
	\caption{Optimized flight trajectory at low wind speed}
	\label{fig:traj_low}
\end{figure}

\begin{figure}
	\includegraphics[width=\linewidth,trim={17cm 3cm 17cm 16cm},clip]{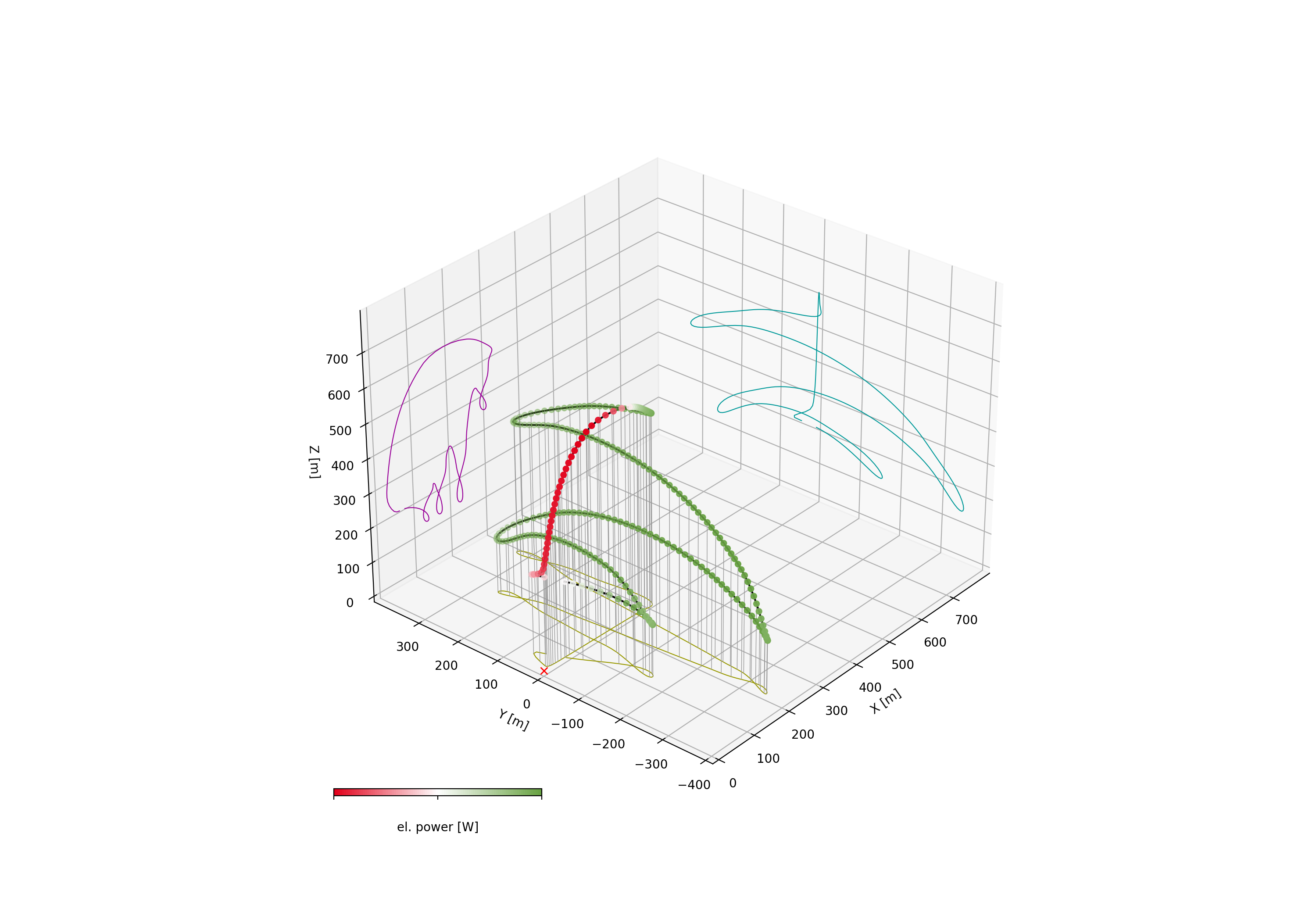}
	\caption{Optimized flight trajectory at high wind speed}
	\label{fig:traj_high}
\end{figure}

The optimization process is fast and delivers high quality, robust results covering  full range of operating conditions, even though the form of the trajectories at low wind speeds, see Fig. \ref{fig:traj_low}, differs significantly from that at high wind speeds, see Fig. \ref{fig:traj_high}. 
At low wind speed, the tether reel-out speed is relatively small. 
Thus the duration of the reel-out phase is high resulting in large number of the well known "figures of eights" flown in this phase, marked using green dots in the figure.
At high wind speeds however, some important constraints limit the lowest kite position and the tether is reeled out much faster thus significantly less time and space is available for the power generation. 
As a result, the corresponding optimal trajectory looks differently.
Nevertheless, it results in significantly higher mean cycle power than that at the low wind speed.    
Importantly, both trajectories are characterized by  short reel-in phase, marked using red dots in the figure, being a significant feature of high efficiency power cycles realized using variable trim kites.

The results of the numerical optimization had a large impact on the R\&D of the PN-14.
On the supervisory control level, the optimal trajectories have been  extensively used to provide very accurate time-varying set points for various feedback control loops of the system.
On the system engineering level, due to the accessibility of the code, its robustness and small computational time of the optimization process, it has been widely used for future system performance projections and various detailed parametric and sensitivity analyses supporting component choice and sizing.




%

%

%

\end{document}